\documentclass[12pt]{article}
\usepackage{color}
\definecolor{darkblue}{rgb}{0.00,0.25,0.50}
\usepackage[colorlinks,filecolor=blue,citecolor=darkblue]{hyperref}

\usepackage{amsfonts,amssymb,amsmath,amsthm}
\usepackage{url}
\usepackage{enumerate}
\usepackage[ukrainian, russian, english]{babel}
\usepackage[cp1251]{inputenc}
\begin{document}\selectlanguage{english}
\thispagestyle{empty}

\title{}

\begin{center}
\textbf{\Large Estimates for approximations by Fourier sums, 
best approximations and best orthogonal  trigonometric approximations  of the classes 
of $(\psi, \beta)$--differentiable functions}
\end{center}
\vskip0.5cm
\begin{center}
A.~S.~Serdyuk${}^1$, T.~A.~Stepanyuk${}^2$\\ \emph{\small
${}^1$Institute of Mathematics NAS of
Ukraine, Kyiv\\
${}^2$Lesya Ukrainka Eastern European National University, Lutsk\\}
\end{center}
\vskip0.5cm

%\address{Institute of Mathematics of the National Academy of
%Sciences of Ukraine\\ 3\\ Tereshenkivska st.\\ 01601\\ Kiev, Ukraine}

\begin{abstract}

We obtain the exact-order estimates for approximations by Fourier sums, best approximations and best orthogonal trigonometric approximations in metrics of spaces   $L_s$, $1\leq s<\infty$,
of classes of
$2\pi$--periodic functions, whose
$(\psi,\beta)$--derivatives belong to unit ball of the space   $L_{\infty}$.
\end{abstract}

\vskip 1.5cm

%%%%%%%%%%%%%%%%%%%%%%%%%%%%%%%%%%%%%%%%%%%%%%%%%%%%%%%%%%%%%%%%%%%%%%%%%

We denote by $L_{p}$,
$1\leq p<\infty$, the space of $2\pi$--periodic functions $f:\mathbb{R}\rightarrow\mathbb{C}$, summable to the power $p$ on  $[0,2\pi)$, in which
the norm is given by the formula
$\|f\|_{p}=\Big(\int\limits_{0}^{2\pi}|f(t)|^{p}dt\Big)^{\frac{1}{p}}$; and we denote by
$L_{\infty}$ the space of
$2\pi$---periodic measurable and essentially bounded functions  $f:\mathbb{R}\rightarrow\mathbb{C}$ with the norm
$\|f\|_{\infty}=\mathop{\rm{ess}\sup}\limits_{t}|f(t)|$;

Let $f:\mathbb{R}\rightarrow\mathbb{R}$ be the function from  $L_{1}$, whose Fourier series has the form
$$
\sum_{k=-\infty}^{\infty}\hat{f}(k)e^{ikx},
$$
where  $\hat{f}(k)=\frac{1}{2\pi}\int\limits_{-\pi}^{\pi}f(t)e^{-ikt}dt$ are Fourier coefficients of the function  $f$,
 $\psi(k)$ is an arbitrary fixed sequence of real numbers and $\beta$
is a fixed real number. Then, if the series
$$
\sum_{k\in \mathbb{Z}\setminus\{0\}}\frac{\hat{f}(k)}{\psi(|k|)}e^{i(kx+\frac{\beta\pi}{2} \mathrm{sign} k)}
$$
\noindent is the Fourier series of some function $\varphi$ from $L_{1}$, then this function is called the $(\psi,\beta)$--derivative of the function $f$  and denoted by
$f_{\beta}^{\psi}$.
A set of functions $f$, whose $(\psi,\beta)$--derivatives exist is denoted by
 $L_{\beta}^{\psi}$  (see \cite{Stepanets}).

If $f\in L^{\psi}_{\beta}$ and, at the same time, $f^{\psi}_{\beta}\in
\mathfrak{N}$, where $\mathfrak{N}\subseteq L_{1}$, then we say that the function $f$ belongs to the class $L^{\psi}_{\beta}\mathfrak{N}$. By $B_{R,p}$ we denote the balls of the radius  $R$ of real--valued functions from $L_p$, i.e., the sets
 $$
 {B_{R,p}:=\left\{\varphi: \ \mathbb{R}\rightarrow\mathbb{R}, \ ||\varphi||_{p}\leq R \right\}}, \ \ R>0, \ \ 1\leq p\leq\infty.
$$

In present paper as $\mathfrak{N}$ we take the unit balls  $B_{1,p}$. Herewith, the functional classes  $L^{\psi}_{\beta}B_{1,p}$
are denoted by
  $L^{\psi}_{\beta,p}$.

In the case $\psi(k)=k^{-r},\ r>0$, the classes  $L_{\beta,p}^{\psi}$ are well--known
Weyl--Nagy classes $W_{\beta,p}^r$.

For functions  $f$ from classes $L^{\psi}_{\beta,p}$ we consider: $L_{s}$--norms of deviations of the functions  $f$ from their partial Fourier sums of order $n-1$, i.e., the quantities
\begin{equation}\label{rho}
\|\rho_{n}(f;\cdot)\|_{s}=
\|f(\cdot)-S_{n-1}(f;\cdot)\|_{s}, \ \ 1\leq s\leq\infty,
\end{equation}
where
$$
S_{n-1}(f;x)=\sum\limits_{k=-n+1}^{n-1}\hat{f}(k)e^{ikx};
$$
best orthogonal trigonometric approximations of the functions $f$ in metric of space  $L_{s}$, i.e., the quantities of the form
\begin{equation}\label{nn_term}
 e^{\bot}_{m}(f)_{s}=
  \inf\limits_{\gamma_{m}}\|f(\cdot)-S_{\gamma_{m}}(f;\cdot) \|_{s},  \ \  \ 1\leq s\leq\infty,
\end{equation}
where $\gamma_{m}$, $m\in\mathbb{N}$, is an arbitrary collection of  $m$ integer numbers, and
$$
S_{\gamma_{m}}(f;x)=\sum\limits_{k\in \gamma_{m}}\hat{f}(k)e^{ikx};
$$
and best approximations of the functions  $f$  in space $L_{s}$, i.e., the quantities of the form
\begin{equation}\label{bestTrig}
  E_{n}(f)_{s}=
\inf\limits_{t_{n-1}\in\mathcal{T}_{2n-1}}\|f-t_{n-1}\|_{s}, \ \
\ 1\leq s\leq \infty,
\end{equation}
where $\mathcal{T}_{2n-1}$ is the subspace of all trigonometric polynomials $t_{n-1}$ with real coefficients  of degrees not greater than ${n-1}$.

We set
\begin{equation}\label{Fsum}
 {\cal E}_{n}(L^{\psi}_{\beta,p})_{s}=
\sup\limits_{f\in
L^{\psi}_{\beta,p}}\|\rho_{n}(f;\cdot)\|_{s}, \ \ 1\leq p,s\leq\infty,
\end{equation}
\begin{equation}\label{n_term}
  e^{\bot}_{n}(L_{\beta,p}^{\psi})_{s}=\sup\limits_{f\in L_{\beta,p}^{\psi}}e^{\bot}_{n}(f)_{s}, \ 1\leq p,s\leq\infty,
\end{equation}
\begin{equation}\label{bestApr}
  {E}_{n}(L^{\psi}_{\beta,p})_{s}=
\sup\limits_{f\in
L^{\psi}_{\beta,p}}E_{n}(f)_{s}, \ \
\ 1\leq p,s\leq \infty.
\end{equation}

The following inequalities follow from given above definitions  (\ref{Fsum})--(\ref{bestApr})
\begin{equation}\label{ineq_comp}
  E_{n}(L_{\beta,p}^{\psi})_{s}\leq {\cal E}_{n}(L^{\psi}_{\beta,p})_{s}, \ \ \
  e^{\bot}_{2n-1}(L_{\beta,p}^{\psi})_{s}\leq {\cal E}_{n}(L^{\psi}_{\beta,p})_{s}, \ 1\leq p,s\leq\infty.
\end{equation}

In present paper we solve the problem about finding the exact order estimates of the quantities  ${\cal E}_{n}(L^{\psi}_{\beta,\infty})_{s}$, ${ E}_{n}(L^{\psi}_{\beta,\infty})_{s}$ and $e^{\perp}_{n}(L^{\psi}_{\beta,\infty})_{s}$ for $1\leq s<\infty$, $\beta\in\mathbb{R}$.

For the Weyl--Nagy classes the exact order estimates of the quantities
 ${\cal E}_{n}(W^{r}_{\beta,p})_{{s}}$ and
${E}_{n}(W^{r}_{\beta,p})_{{s}}$ are known  for all admissible values of parameters   $r$, $p$, $s$ and  $\beta$, i.e., for $r>\max\big\{\frac{1}{p}-\frac{1}{s}, 0\big\}$,
$\beta\in\mathbb{R}$ and ${1\leq p,s\leq\infty}$ (see, e.g., \cite[p. 47--49]{T}).
What concerning the  best orthogonal trigonometric approximations    $e^{\bot}_{n}(W^{r}_{\beta,p})_s$, so order estimates are known for them  (see \cite{Belinsky1988}--\cite{Romanyuk2012})
for various   (but not for all possible) values of the parameters  $r,p,s$ and $\beta$.

Order estimates of the quantities  (\ref{Fsum})--(\ref{bestApr}) under certain restrictions for the parameters  $r,p,s$ and $\beta$  were established in the works  \cite{Stepanets}, \cite{Kuchpel}--\cite{Serdyuk_Stepaniuk10}. However, the case  $p=\infty$, ${1\leq s\leq\infty}$ for some or another reasons hasn't been investigated yet.

We denote by $P$ the set of positive, almost decreasing sequences  $\psi(k)$, $k\geq1$, (we remind, that sequence $\psi(k)$ almost decreases, if there exists a positive constant  $M$ such that for arbitrary
 $k_{1}\leq k_{2}$ the following inequality is satisfied
$\psi(k_{2})\leq M\psi(k_{1})$)
such that
$$
\sup\limits_{m\in\mathbb{N}}\sum\limits_{k=2^{m}}^{2^{m+1}}|\psi_n(k+1)-\psi_n(k)|\leq K\psi(n),
$$
where
$$
\psi_n(k)={\left\{\begin{array}{cc}
0, \ \ \  \ & k<n, \\
\psi(k),   & k\geq n, \
  \end{array} \right.}
$$
and $K$ is the quantity uniformly bounded with respect to  $n$.

\textbf{Theorem 1.} \emph{ Let  $\psi\in P$, ${1\leq s<\infty}$ and
$\beta\in \mathbb{R}$. Then}
\begin{equation}\label{theorem1}
 E_{n}(L^{\psi}_{\beta,\infty})_{s}\asymp \mathcal{E}_{n}(L^{\psi}_{\beta,\infty})_{s}\asymp
\psi(n).
\end{equation}
Here and in what follows, we write  ${A(n)\asymp B(n)}$ for postive sequences  $A(n)$ and $B(n)$ to denote that there are positive
constants $K_{1}$ and
$K_{2}$ such that ${K_{1}B(n)\leq A(n)\leq K_{2}B(n)}$,
$n\in\mathbb{N}$.

\emph{Proof.} At first let's  prove that the following inequality is true
\begin{equation}\label{q4}
\mathcal{E}_{n}(L^{\psi}_{\beta,\infty})_{s}\leq K^{(1)}
\psi(n), \ \ 1\leq s<\infty.
\end{equation}
In inequality (\ref{q4}) and henceforth by
 $K^{(i)}$, $i=1,2,...$ we denote quantities uniformly bounded with respect to   $n$.

If $f\in  L^{\psi}_{\beta,\infty}$, then
\begin{equation}\label{st}
  \|f^{\psi}_{\beta}\|_{s}\leq (2\pi)^{\frac{1}{s}}\|f^{\psi}_{\beta}\|_{\infty}\leq (2\pi)^{\frac{1}{s}},
\end{equation}
and so, it is obviously that
\begin{equation}\label{00q100}
 L^{\psi}_{\beta,\infty}\subset L^{\psi}_{\beta}B_{(2\pi)^{\frac{1}{s}}, s}\subset L^{\psi}_{\beta}L_{s},
 \ 1\leq s< \infty.
\end{equation}

The following proposition follows from the theorem  6.7.1 in \cite{Stepanets}.

\textbf{\emph{Proposition 1.}} {\it Let $1<s<\infty$,  $\psi\in P$, $f\in L^{\psi}_{\beta} L_{s}$ and  $\beta\in\mathbb{R}$. Then for arbitrary  $n\in\mathbb{N}$
 there exists a positive constant  $K$, which is uniformly bounded with respect to  $n$ and $f$ and such that }
\begin{equation}\label{statement1}
 \|\rho_{n}(f;x)\|_s\leq K\psi(n)E_{n}(f^{\psi}_{\beta})_{s}.
\end{equation}

Taking into account (\ref{st}), (\ref{00q100}) and in view of proposition 1, we obtain the following estimates
\begin{equation}\label{qq01}
  \mathcal{E}_{n}(L^{\psi}_{\beta,\infty})_{s}\leq
   \mathcal{E}_{n}\big(L^{\psi}_{\beta}B_{(2\pi)^{\frac{1}{s}},s}\big)_{s}\leq
 (2\pi)^{\frac{1}{s}}K\psi(n), \ \ 1<s<\infty.
\end{equation}

Thus, the inequalities  (\ref{q4})   are proved for $1<s<\infty$.

Let's show the rightness of correlation  (\ref{q4}) for $s=1$.
We use the following statement (see, e.g., \cite[p. 8]{T}).

\textbf{Proposition 2.} \emph{Let $1\leq q\leq p\leq \infty$. On this if $f\in L_p$, then $f\in L_q$ and}
\begin{equation}\label{statement2}
\|f\|_{q}\leq \left(2\pi\right)^{\frac{1}{q}-\frac{1}{p}}\|f\|_{p}.
\end{equation}

By using (\ref{statement2}) for $q=1$, $p=2$ and inequality (\ref{qq01}) for $s=2$, we obtain
$$
 \mathcal{E}_{n}(L^{\psi}_{\beta,\infty})_{1}=\sup\limits_{f\in
L^{\psi}_{\beta,\infty}}\|f(\cdot)-S_{n-1}(f;\cdot)\|_{1}\leq
$$
\begin{equation}\label{q1}
 \leq (2\pi)^{\frac{1}{2}}\sup\limits_{f\in
L^{\psi}_{\beta,\infty}}\|f(\cdot)-S_{n-1}(f;\cdot)\|_{2}=(2\pi)^{\frac{1}{2}} \mathcal{E}_{n}(L^{\psi}_{\beta,\infty})_{2}\leq K^{(1)}\psi(n).
\end{equation}
The rightness of the inequality (\ref{q4}) follows from  (\ref{qq01}) and  (\ref{q1}).

 To obtain the lower bound of the quantity  $E_{n}(L^{\psi}_{\beta,\infty})_{s}$, we consider the following function
$$
f_{1}(t)=f_1(\psi;n;t)=\psi(n)\cos nt.
$$
It is obviously, that  $f_{1}\in L^{\psi}_{\beta,\infty}$ and
$f_{1}\perp t_{n-1}$ for arbitrary
 $t_{n-1}\in~ \mathcal{T}_{2n-1}$.
 Therefore
 \begin{equation}\label{f3}
\int\limits_{-\pi}^{\pi}(f_{1}(t)-t_{n-1}(t))\cos
ntdt=\int\limits_{-\pi}^{\pi}f_{1}(t)\cos ntdt=\pi\psi(n)  \ \  \forall t_{n-1}\in~ \mathcal{T}_{2n-1}.
 \end{equation}

On the other hand,  taking into account the proposition 2 for $q=1$, $p=s$, we get
$$
\int\limits_{-\pi}^{\pi}(f_{1}(t)-t_{n-1}(t))\cos ntdt\leq \|f_{1}-t_{n-1}\|_{1}\leq
$$
\begin{equation}\label{f4}
\leq (2\pi)^{1-\frac{1}{s}}\|f_{1}-t_{n-1}\|_{s}, \ \ 1\leq s\leq\infty, \ \ \forall t_{n-1}\in \mathcal{T}_{2n-1}.
\end{equation}

In view of (\ref{f3})--(\ref{f4}) we arrive at the inequalities
\begin{equation}\label{q5}
 E_{n}(L^{\psi}_{\beta,\infty})_{s}\geq E_{n}(f_{1})_{s}=\inf\limits_{t_{n-1}\in \mathcal{T}_{2n-1}}\|f_{1}-t_{n-1}\|_{s}\geq\frac{1}{2}\psi(n), \ \ 1\leq s\leq\infty.
\end{equation}
Theorem 1 is proved.

We denote by  $B$ the set of positive sequences
  $\psi(k)$, $k\in\mathbb{N}$,  for each of which there exists a positive constant $K$ such that $
{\frac{\psi(k)}{\psi(2k)}\leq K}$, $k\in\mathbb{N} $.
The sequences ${\psi(k)=k^{-r}}$, ${r>0}$,
 $\psi(k)=\ln^{-\varepsilon}(k+1), \ \varepsilon>0$, etc. are representatives of the set
 $B$.

\textbf{Theorem 2.} \emph{ Let $\psi\in P\cap B$, ${1\leq s<\infty}$ and
$\beta\in \mathbb{R}$. Then}
\begin{equation}\label{theorem2}
  e_{2n}^{\perp}(L^{\psi}_{\beta,\infty})_{s}\asymp e_{2n-1}^{\perp}(L^{\psi}_{\beta,\infty})_{s}\asymp
\psi(n).
\end{equation}
\emph{Proof.} It is follows from the formulas  (\ref{ineq_comp}) and (\ref{q4}), that under the conditions of the theorem 1, next inequalities are true
\begin{equation}\label{stst}
  e_{2n}^{\perp}(L^{\psi}_{\beta,\infty})_{s}\leq  e_{2n-1}^{\perp}(L^{\psi}_{\beta,\infty})_{s}\leq  \mathcal{E}_{n}(L^{\psi}_{\beta,\infty})_{s}\leq K^{(1)}\psi(n).
\end{equation}
Now we determine a lower bound of the quantity  $e_{2n}^{\perp}(L^{\psi}_{\beta,\infty})_{s}$.
For this we use the well--known result of Rudin--Shapiro  (see, e.g., lemma 6.32.1 in \cite{Djachenko_Ulianov}).

\textbf{Proposition 3.} \emph{There exists sequence of numbers  $\{\varepsilon_{k}\}_{k=0}^{\infty}$, such that $\varepsilon_{k}=\pm1$ and}
\begin{equation}\label{statement5}
\Big\|\sum\limits_{k=0}^{m}\varepsilon_{k}e^{ikx}\Big\|_{\infty}\leq 5\sqrt{m+1}, \ \ \ m=0,1,...
\end{equation}

Taking into account proposition 3 for  $m=2n-1$, we choose the sequence of numbers  $\{\xi_{k}\}_{k=0}^{\infty}$, $\xi_{k}=\pm1$ such that
\begin{equation}\label{q8}
 \Big\|\sum\limits_{k=0}^{2n-1}\xi_{k}e^{ikx}\Big\|_{\infty}\leq 5\sqrt{2n}.
\end{equation}

We set
$$
\psi(0):=\psi(1)
$$
and consider the function
\begin{equation}\label{q6}
f_{2}(t)=f_{2}(\psi;n;t):=
\frac{1}{10\sqrt{2n}+2}\sum\limits_{k=-2n+1}^{2n-1}\xi_{|k|}\psi(|k|)e^{ikt}.
\end{equation}

Since, according to definition of  $(\psi,\beta)$--derivative and the inequality  (\ref{q8}),
$$
\big\|(f_{2})^{\psi}_{\beta}\big\|_{\infty}=
\frac{1}{10\sqrt{2n}+2}\Big\|\sum\limits_{k=1}^{2n-1}\xi_{k}e^{i(kt+\frac{\beta\pi}{2})}+ \sum\limits_{k=1}^{2n-1}\xi_{k}e^{i(-kt-\frac{\beta\pi}{2})}\Big\|_{\infty}\leq
$$
$$
\leq \frac{1}{10\sqrt{2n}+2}\Big(\Big\|\sum\limits_{k=1}^{2n-1}\xi_{k}e^{i(kt+\frac{\beta\pi}{2})}\Big\|_{\infty}+ \Big\|\sum\limits_{k=1}^{2n-1}\xi_{k}e^{i(-kt-\frac{\beta\pi}{2})}\Big\|_{\infty}\Big)=
$$
$$
=\frac{1}{5\sqrt{2n}+1}\Big\|\sum\limits_{k=1}^{2n-1}\xi_{k}e^{ikt} \Big\|_{\infty}
\leq 1,
$$
so $f_{2}\in L^{\psi}_{\beta,\infty}$.

We consider the quantity
$$
I=\inf\limits_{\gamma_{2n}}\bigg|\int\limits_{-\pi}^{\pi}(f_{2}(t)-S_{\gamma_{2n}}(f_{2};t))
\sum\limits_{k=-2n+1}^{2n-1}\xi_{|k|}e^{ikt}dt\bigg|.
$$
 By virtue of Holder's inequality, proposition  2 and correlation  (\ref{q8}) for $1\leq s<\infty$, $\frac{1}{s}+\frac{1}{s'}=1$
 $$
 I\leq\inf\limits_{\gamma_{2n}}\|f_{2}(t)-S_{\gamma_{2n}}(f_{2};t)\|_{s}\Big\|\sum\limits_{k=-2n+1}^{2n-1}\xi_{|k|}e^{ikt}\Big\|_{s'}=
 $$
 $$
=
 e^{\bot}_{2n}(f_{2})_{s}\Big\|\sum\limits_{k=-2n+1}^{2n-1}\xi_{|k|}e^{ikt}\Big\|_{s'}
\leq (2\pi)^{\frac{1}{s'}}e^{\bot}_{2n}(f_{2})_{s}\Big\|\sum\limits_{k=-2n+1}^{2n-1}\xi_{|k|}e^{ikt}\Big\|_{\infty}\leq
 $$
 $$
 \leq2\pi
e^{\bot}_{2n}(f_{2})_{s}\Big(\Big\|\sum\limits_{k=0}^{2n-1}\xi_{k}e^{ikt}\Big\|_{\infty}+\Big\|\sum\limits_{k=1}^{2n-1}\xi_{k}e^{-ikt}\Big\|_{\infty}\Big)
\leq
 $$
\begin{equation}\label{q9}
\leq2\pi
e^{\bot}_{2n}(f_{2})_{s}\Big(2\Big\|\sum\limits_{k=0}^{2n-1}\xi_{k}e^{ikt}\Big\|_{\infty}+1\Big)
\leq 2\pi(10\sqrt{2n}+1) \ e^{\bot}_{2n}(f_{2})_{s}.
\end{equation}

On the other hand, taking into account the orthogonality of trigonometric system  $\{e^{ikt}\}$ and the fact that $\xi_{k}^{2}=1$, we obtain
$$
I=\frac{1}{10\sqrt{2n}+2}\inf\limits_{\gamma_{2n}}\bigg|\int\limits_{-\pi}^{\pi}
{\mathop{\sum}\limits_{ |k|\leq 2n-1,\atop k \notin\gamma_{2n} }}\xi_{|k|}\psi(|k|)e^{ikt}\sum\limits_{k=-2n+1}^{2n-1}\xi_{|k|}e^{ikt}dt\bigg|=
$$
\begin{equation}\label{q10}
=\frac{\pi}{5\sqrt{2n}+1}\inf\limits_{\gamma_{2n}}
{\mathop{\sum}\limits_{ |k|\leq 2n-1,\atop k \notin\gamma_{2n} }}\psi(|k|).
\end{equation}

Since the sequence  $\psi(k)$ almost decreases, so
\begin{equation}\label{pp}
 \inf\limits_{\gamma_{2n}}
{\mathop{\sum}\limits_{ |k|\leq 2n-1,\atop k \notin\gamma_{2n} }}\psi(|k|)
\geq K^{(2)}\inf\limits_{\gamma_{2n}}
{\mathop{\sum}\limits_{ |k|\leq 2n-1,\atop k \notin\gamma_{2n} }}\psi(2n-1)=K^{(2)}\psi(2n-1)(2n-1).
\end{equation}

In view of (\ref{q9})--(\ref{pp})  we get
\begin{equation}\label{p}
  e^{\bot}_{2n}(f_{2})_{s}\geq\frac{K^{(2)}\psi(2n-1)(2n-1)}{(10\sqrt{2n}+2)(10\sqrt{2n}+1)}\geq K^{(3)}\psi(2n).
  \end{equation}

Since, if $\psi\in B$, so $\psi(2n)\geq K^{(4)}\psi(n)$, and, hence, taking into account  (\ref{p}), we find
\begin{equation}\label{q11}
  e_{2n}^{\perp}(L^{\psi}_{\beta,\infty})_{s}\geq e^{\bot}_{2n}(f_{2})_{s}\geq K^{(5)}\psi(n).
\end{equation}

Estimates (\ref{theorem2}) follow from
 (\ref{stst}) and (\ref{q11}).
Theorem 2 is proved.

\textbf{Corollary 1.} \emph{ Let $r>0$, ${1\leq s<\infty}$ and
$\beta\in \mathbb{R}$. Then}
\begin{equation}\label{conseq1}
  e_{2n}^{\perp}(W^{r}_{\beta,\infty})_{s}\asymp e_{2n-1}^{\perp}(W^{r}_{\beta,\infty})_{s}\asymp
n^{-r}.
\end{equation}

E-mail: \href{mailto:serdyuk@imath.kiev.ua}{serdyuk@imath.kiev.ua},
\href{mailto:tania_stepaniuk@ukr.net}{tania$_{-}$stepaniuk@ukr.net}


\newpage
\begin{thebibliography}{10}

\bibitem{Stepanets}
{\sc A.I. Stepanets,}
\emph{Methods of Approximation Theory}, VSP: Leiden, Boston,  (2005).



\bibitem{T}
 {\sc V.N. Temlyakov,} \emph{Approximation of Periodic Function}: Nova
Science Publichers, Inc., (1993).




\bibitem{Belinsky1988}
{\sc E.S. Belinsky},
\emph{Approximation by a "floating" system of exponents on the classes of periodic functions with bounded mixed derivative}, Issled. po teorii func. mnog. vesch. perem., Jaroslavl':  Jaroslav. un--t.,
(1988), 16--33.  [in Russian]

\bibitem{Temljakov1986}
{\sc V.N. Temlyakov,} \emph{Approximations of functions with bounded mixed derivative}, Trudy Mat. Inst. Steklov., \textbf{178} (1986), 3–113 [in Russian]; \textbf{English translation:}
Proceedings of the Steklov Institute of Mathematics, \textbf{178} (1989), 1–121.	



\bibitem{Kashyn_Teml}
{\sc B.S. Kashin and V.N. Temlyakov,} \emph{On best $m$--term approximations and the entropy of sets in the space $L_1$}, Mat. Zametki, \textbf{56}:5 (1994), 57--86 [in Russian];
 \textbf{English translation:} Mathematical Notes, \textbf{56}:5 (1994), 1137–1157.


\bibitem{Romanyuk1992_4}
{\sc A.S. Romanyuk,}
\emph{Approximation of the classes $B_{p,\Theta}^{r}$ of periodic functions of several variables by partial Fourier sums with arbitrary choice of harmonics}, Proc. of the Institute of Mathematics "Fourier series: theory and applications", Ukrainian National Academy of Sciences [in Russian],
Kyiv,  (1992),  112--118.

\bibitem{Romanyuk2002}
{\sc A.S. Romanyuk,} \emph{Approximation of classes of periodic functions in several variables}, Mat. Zametki, \textbf{71}:1 (2002), 109–121 [in Russian]; \textbf{English translation:}
Mathematical Notes, \textbf{71}:1 (2002), 98–109


\bibitem{Romanyuk2007}
{\sc A.S. Romanyuk,} \emph{Best trigonometric approximations for some classes of periodic functions of several variables in the uniform metric},
 Mat. Zametki, \textbf{82}:2 (2007), 247–261 [in Russian]; \textbf{English translation:} Mathematical Notes, {\bf 82}:2 (2007), 216--228.


\bibitem{Romanyuk2012}
 {\sc A.S. Romanyuk,}
\emph{Approximation characteristics of classes of periodic functions of several variables}, Proc. of the Institute of Mathematics, Ukrainian National Academy of Sciences [in Russian],
Kyiv,  \textbf{93} (2012). [in Russian]



\bibitem{Kuchpel}
{\sc A.I. Stepanets and A.K. Kushpel’,} \emph{Convergent rate of Fourier series and best approximations in the space $L_{p}$},  Ukr. Mat. Zh., \textbf{39}:4 (1987), 483--492 [in Russian]; \textbf{English translation:} Ukr. Math. J., \textbf{39}:4 (1987), 389–398.





\bibitem{Serdyuk_grabova}
{\sc U.Z. Hrabova and A.S. Serdyuk,} \emph{Order estimates for the best approximations and approximations by Fourier sums of the classes of
$(\psi,\beta)$--differential functions}, Ukr. Mat. Zh., \textbf{65}:9, 1186--1197 (2013) [in Ukrainian]; \textbf{English translation:} Ukr. Math. J., \textbf{65}:9 (2013), 1319--1331.



\bibitem{S_S}
{\sc A.S. Serdyuk and T.A. Stepanyuk,} \emph{Order estimates for the best approximations and approximations of the classes of infinitely differential functions by Fourier sums}, Proc. of the Institute of Mathematics, Ukrainian National Academy of Sciences [in Ukrainian],
Kyiv, \textbf{10}:1 (2013),  255–282.

\bibitem{Serdyuk_Stepaniuk7}
{\sc A.S. Serdyuk and T.A. Stepanyuk,} Estimates for the best approximations of the classes of infinitely differentiable functions in
uniform and integral metrics, Ukr. Mat. Zh., \textbf{66}:9 (2014), 1244–1256 [in Ukrainian]; \textbf{English translation:} Ukr. Math. J., \textbf{66}:9 (2015), 1393--1407.


\bibitem{Stepaniuk2014}
{\sc T.A. Stepanyuk,} \emph{Estimates for the best approximations and approximations by Fourier sums of the classes of  convolutions of periodic
functions of not high smoothness in the integral metrics}, Proc. of the Institute of Mathematics, Ukrainian National Academy of Sciences [in Ukrainian],
Kyiv, \textbf{11}:3 (2014),  241–269.


\bibitem{Serdyuk_Stepaniuk8}
{\sc A.S. Serdyuk and T.A. Stepanyuk,}  \emph{Order estimates for the best approximations and approximations
by Fourier sums in the classes of convolutions of periodic
functions of low smoothness in the uniform metric}, Ukr. Mat. Zh., \textbf{66}:12 (2014), 1658--1675 [in Ukrainian]; \textbf{English translation:} Ukr. Math. J., \textbf{66}:12 (2015), 1862--1882.





\bibitem{Fedorenko1999}
{\sc O.S. Fedorenko,} \emph{On the best  $m$--term trigonometric  and orthogonal  trigonometric  approximations}, Ukr. Mat. Zh., \textbf{51}:12 (1999), 1719--1721 [in Ukrainian]; \textbf{English translation:} Ukr. Math. J., \textbf{51}:12 (2000), 1945--1949.



\bibitem{Fedorenko2000}
{\sc O.S. Fedorenko,} \emph{Approximation of
 $(\psi,\beta)$--differentiable functions by trigonometric polynomials}: Autoref. dissertation ... cand. phys.--math.  sciences, Kyiv: Insitute of Mathematics of NAS of Ukraine (2001). [in Ukrainian]

\bibitem{Shkapa}
{\sc V.V. Shkapa,}
 \emph{Best  orthogonal trigonometric approximations of functions of classes $L^{\psi}_{\beta,1}$},
Proc. of the Institute of Mathematics, Ukrainian National Academy of Sciences [in Ukrainian],
Kyiv, \textbf{11}:3 (2014),  315--329.


\bibitem{Shkapa1}
{\sc V.V. Shkapa,}
\emph{Estimates of the best $M$--term and orthogonal trigonometric approximations of functions of classes $L^{\psi}_{\beta,p}$  in the uniform metrics},
Proc. of the Institute of Mathematics, Ukrainian National Academy of Sciences [in Ukrainian],
Kyiv, \textbf{11}:2 (2014),  305--307.



\bibitem{Serdyuk_Stepaniuk10}
{\sc A.S. Serdyuk and T.A. Stepanyuk,}  O\emph{rder estimates for the best orthogonal trigonometric approximations of the classes of convolutions of periodic
functions of low smoothness in the uniform metric}, Ukr. Mat. Zh., \textbf{66}:12 (2014), 1658--1675 [in Ukrainian].


\bibitem{Djachenko_Ulianov}
{\sc M.I. Dyachenko and P.L. Ulyanov,} \emph{Measure and integral}, "Factorial",  Moscow (1998). [in Russian]
\end{thebibliography}
\end{document}